# HIGH-DIMENSIONAL GENERALIZED LINEAR MODELS AND THE LASSO

By Sara A. van de Geer

*ETH Zürich*

We consider high-dimensional generalized linear models with Lipschitz loss functions, and prove a nonasymptotic oracle inequality for the empirical risk minimizer with Lasso penalty. The penalty is based on the coefficients in the linear predictor, after normalization with the empirical norm. The examples include logistic regression, density estimation and classification with hinge loss. Least squares regression is also discussed.

**1. Introduction.** We consider the lasso penalty for high-dimensional generalized linear models. Let $Y \in \mathcal{Y} \subset \mathbf{R}$ be a real-valued (response) variable and $X$ be a co-variable with values in some space $\mathcal{X}$. Let

$$\mathcal{F} = \left\{ f_\theta(\cdot) = \sum_{k=1}^m \theta_k \psi_k(\cdot), \theta \in \Theta \right\}$$

be a (subset of a) linear space of functions on $\mathcal{X}$. We let $\Theta$ be a convex subset of $\mathbf{R}^m$, possibly $\Theta = \mathbf{R}^m$. The functions $\{\psi_k\}_{k=1}^m$ form a given system of real-valued base functions on $\mathcal{X}$.

Let $\gamma_f : \mathcal{X} \times \mathcal{Y} \to \mathbf{R}$ be some loss function, and let $\{(X_i, Y_i)\}_{i=1}^n$ be i.i.d. copies of $(X, Y)$. We consider the estimator with lasso penalty

$$\hat{\theta}_n := \arg\min_{\theta \in \Theta} \left\{ \frac{1}{n} \sum_{i=1}^n \gamma_{f_\theta}(X_i, Y_i) + \lambda_n \hat{I}(\theta) \right\},$$

where

$$\hat{I}(\theta) := \sum_{k=1}^m \hat{\sigma}_k |\theta_k|$$

denotes the weighted $\ell_1$ norm of the vector $\theta \in \mathbf{R}^m$, with random weights

$$\hat{\sigma}_k := \left( \frac{1}{n} \sum_{i=1}^{n} \psi_k^2(X_i) \right)^{1/2}.$$

Moreover, the smoothing parameter $\lambda_n$ controls the amount of complexity regularization.

Indeed, when $m$ is large (possibly $m > n$), a complexity penalty is needed to avoid overfitting. One could penalize the number of nonzero coefficients, that is, use an $\ell_0$ penalty, but this leads to a nonconvex optimization problem. The lasso is a convex $\ell_1$ penalty which often behaves similarly as the $\ell_0$ penalty. It corresponds to soft thresholding in the case of quadratic loss and orthogonal design, see Donoho (1995). In Donoho (2006a), (2006b), the agreement of $\ell_1$- and $\ell_0$-solutions in general (nonorthogonal) systems is developed further. We will address the distributional properties of the lasso penalty. The acronym "lasso" (least absolute shrinkage and selection operator) was introduced by Tibshirani (1996), with further work in Hastie, Tibshirani and Friedman (2001).

Let $P$ be the distribution of $(X, Y)$. The target function $\bar{f}$ is defined as

$$\bar{f} := \underset{f \in \mathbf{F}}{\arg \min}\, P\gamma_f,$$

where $\mathbf{F} \supseteq \mathcal{F}$ (and assuming for simplicity that there is a unique minimum). We will show that if the target $\bar{f}$ can be well approximated by a sparse function $f_{\theta_n^*}$, that is a function $f_{\theta_n^*}$ with only a few nonzero coefficients $\theta_{n,k}^*$, the estimator $\hat{\theta}_n$ will have prediction error roughly as if it knew this sparseness. In this sense, the estimator mimics a sparseness oracle. Our results are of the same spirit as those in Bunea, Tsybakov and Wegkamp (2007b), which we learned about when writing this paper, and which is about quadratic loss. We will assume Lipschitz loss (see Assumption L below). Examples are the loss functions used in quantile regression, logistic regression, the hinge loss function used in classification, and so forth. Quadratic loss can, however, be handled using similar arguments, see Example 4.

Let us briefly sketch our main result. The excess risk of $f$ is

$$\mathcal{E}(f) := P\gamma_f - P\gamma_{\bar{f}}.$$

We will derive a probability inequality for the excess risk $\mathcal{E}(f_{\hat{\theta}_n})$ (see Theorems 2.1 and 2.2). This inequality implies that with large probability,

$$\mathcal{E}(f_{\hat{\theta}_n}) \leq \text{const.} \times \min_{\theta \in \Theta} \{ \mathcal{E}(f_\theta) + \mathcal{V}_\theta \}.$$

Here, the error term $\mathcal{V}_\theta$ will be referred to as "estimation error." It is typically proportional to $\lambda_n^2 \times \dim_\theta$, where $\dim_\theta := |\{\theta_k \neq 0\}|$. Moreover, as we will see, the smoothing parameter $\lambda_n$ can be chosen of order $\sqrt{\log m/n}$. Thus, again typically, the term $\mathcal{V}_\theta$ is of order $\log m \times \dim_\theta/n$, which is the usual form for the estimation error when estimating $\dim_\theta$ parameters, with the $(\log m)$—term the price to pay for not knowing beforehand which parameters are relevant.

We will study a more general situation, with general margin behavior, that is, the behavior of the excess risk near the target (see Assumption B). The term $\mathcal{V}_\theta$ will depend on this margin behavior as well as on the dependence structure of the features $\{\psi_k\}_{k=1}^m$ (see Assumption C).

The typical situation sketched above corresponds to "quadratic" margin behavior and, for example, a well conditioned inner product matrix of the features.

To avoid digressions, we will not discuss in detail in the case where the inner product matrix is allowed to be singular, but only briefly sketch an approach to handle this (see Section 3.1).

There are quite a few papers on estimation with the lasso in high-dimensional situations. Bunea, Tsybakov and Wegkamp (2006) and Bunea, Tsybakov and Wegkamp (2007b) are for the case of quadratic loss. In Tarigan and van de Geer (2006), the case of hinge loss function is considered, and adaptation to the margin in classification. Greenshtein (2006) studies general loss with $\ell_1$ constraint. The group lasso is studied and applied in Meier, van de Geer and Bühlmann (2008) for logistic loss and Dahinden et al. (2008) for log-linear models. Bunea, Tsybakov and Wegkamp (2007a) studies the lasso for density estimation.

We will mainly focus on the prediction error of the estimator, that is, its excess risk, and not so much on selection of variables. Recent work where the lasso is invoked for variable selection is Meinshausen (2007), Meinshausen and Bühlmann (2006), Zhao and Yu (2006), Zhang and Huang (2006) and Meinshausen and Yu (2007). As we will see, we will obtain inequalities for the $\ell_1$ distance between $\hat{\theta}_n$ and the "oracle" $\theta_n^*$ defined below [see (6) and (9)]. These inequalities can be invoked to prove variable selection, possibly after truncating the estimated coefficients.

We extend results from the papers [Loubes and van de Geer (2002), van de Geer (2003)], in several directions. First, the case of random design is considered. The number of variables $m$ is allowed to be (much) larger than the number of observations $n$. The results are explicit, with constants with nonasymptotic relevance. It is not assumed a priori that the regression functions $f_\theta \in \mathcal{F}$ are uniformly bounded in sup norm. Convergence of the estimated coefficients $\hat{\theta}_n$ in $\ell_1$ norm is obtained. And finally, the penalty is based on the *weighted* $\ell_1$ norm of the coefficients $\theta$, where it is allowed that the base functions $\psi_k$ are normalized with their *empirical* norm $\hat{\sigma}_k$.

The paper is organized as follows. In the remainder of this section, we introduce some notation and assumptions that will be used throughout the paper. In Section 2, we first consider the case where the weights $\sigma_k := (E\psi_k^2(X))^{1/2}$ are known. Then the results are more simple to formulate. It serves as a preparation for the situation of unknown $\sigma_k$. We have chosen to present the main results (Theorems 2.1 and 2.2) with rather arbitrary values for the constants involved, so that their presentation is more transparent. The explicit dependence on the constants can be found in Theorem A.4 (for the case of known $\sigma_k$) and Theorem A.5 (for the case of unknown $\sigma_k$).

Section 3 discusses the Assumptions L, A, B and C given below, and presents some examples. Some extensions are considered, for instance hinge loss, or support vector machine loss, (which may need an adjustment of Assumption B), and quadratic loss (which is not Lipschitz). Moreover, we discuss the case where some coefficients (for instance, the constant term) are not penalized.

The proof of the main theorems is based on the concentration inequality of Bousquet (2002). We moreover use a convexity argument to obtain reasonable constants. All proofs are deferred to the Appendix.

We use the following notation. The empirical distribution based on the sample $\{(X_i, Y_i)\}_{i=1}^n$ is denoted by $P_n$, and the empirical distribution of the covariates $\{X_i\}_{i=1}^n$ is written as $Q_n$. The distribution of $X$ is denoted by $Q$. We let $\sigma_k^2 := Q\psi_k^2$ and $\hat{\sigma}_k^2 := Q_n\psi_k^2$, $k = 1, \ldots, m$. The $L_2(Q)$ norm is written as $\|\cdot\|$. Moreover, $\|\cdot\|_\infty$ denotes the sup norm.

We impose four basic assumptions: Assumptions L, A, B and C.

ASSUMPTION L.   The loss function $\gamma_f$ is of the form $\gamma_f(x,y) = \gamma(f(x), y) + b(f)$, where $b(f)$ is a constant which is convex in $f$, and $\gamma(\cdot, y)$ is convex for all $y \in \mathcal{Y}$. Moreover, it satisfies the Lipschitz property

$$|\gamma(f_\theta(x), y) - \gamma(f_{\tilde{\theta}}(x), y)| \leq |f_\theta(x) - f_{\tilde{\theta}}(x)|$$

$$\forall\, (x,y) \in \mathcal{X} \times \mathcal{Y},\ \forall\, \theta, \tilde{\theta} \in \Theta.$$

Note that by a rescaling argument, there is no loss of generality that Assumption L takes the Lipschitz constant equal to one.

ASSUMPTION A.   It holds that

$$K_m := \max_{1 \leq k \leq m} \frac{\|\psi_k\|_\infty}{\sigma_k} < \infty.$$

ASSUMPTION B.   There exists an $\eta > 0$ and strictly convex increasing $G$, such that for all $\theta \in \Theta$ with $\|f_\theta - \bar{f}\|_\infty \leq \eta$, one has

$$\mathcal{E}(f_\theta) \geq G(\|f_\theta - \bar{f}\|).$$

ASSUMPTION C.   There exists a function $D(\cdot)$ on the subsets of the index set $\{1, \ldots, m\}$, such that for all $\mathcal{K} \subset \{1, \ldots, m\}$, and for all $\theta \in \Theta$ and $\tilde{\theta} \in \Theta$, we have

$$\sum_{k \in \mathcal{K}} \sigma_k |\theta_k - \tilde{\theta}_k| \leq \sqrt{D(\mathcal{K})} \|f_\theta - f_{\tilde{\theta}}\|.$$

The convex conjugate of the function $G$ given in Assumption B is denoted by $H$ [see, e.g., Rockafeller (1970)]. Hence, by definition, for any positive $u$ and $v$, we have

$$uv \leq G(u) + H(v).$$

We define moreover for all $\theta \in \Theta$,

$$D_\theta := D(\{k : |\theta_k| \neq 0\}),$$

with $D(\cdot)$ given in Assumption C.

We let

$$(1) \qquad \bar{a}_n = 4a_n, \qquad a_n := \left( \sqrt{\frac{2 \log(2m)}{n}} + \frac{\log(2m)}{n} K_m \right),$$

with $K_m$ the bound given in Assumption A. We further let for $t > 0$,

$$(2) \qquad \lambda_{n,0} := \lambda_{n,0}(t) := a_n \left( 1 + t \sqrt{2(1 + 2a_n K_m)} + \frac{2t^2 a_n K_m}{3} \right)$$

and

$$(3) \qquad \bar{\lambda}_{n,0} := \bar{\lambda}_{n,0}(t) := \bar{a}_n \left( 1 + t \sqrt{2(1 + 2\bar{a}_n K_m)} + \frac{2t^2 \bar{a}_n K_m}{3} \right).$$

The quantity $\bar{\lambda}_{n,0}$ will be a lower bound for the value of the smoothing parameter $\lambda_n$. Its particular form comes from Bousquet's inequality (see Theorem A.1). The choice of $t$ is "free." As we will see, large values of $t$ may give rise to large excess risk of the estimator, but more "confidence" in the upper bound for the excess risk. In Section 3.2, we will in fact fix the value of $\bar{\lambda}_{n,0}(t)$, and the corresponding value of $t$ may be unknown. The latter occurs when $K_m$ is unknown.

Let

$$I(\theta) := \sum_{k=1}^m \sigma_k |\theta_k|.$$

We call $I(\theta)$ the (theoretical) $\ell_1$ norm of $\theta$, and $\hat{I}(\theta) = \sum_{k=1}^m \hat{\sigma}_k |\theta_k|$ its empirical $\ell_1$ norm. Moreover, for any $\theta$ and $\tilde{\theta}$ in $\Theta$, we let

$$I_1(\theta | \tilde{\theta}) := \sum_{k : \tilde{\theta}_k \neq 0} \sigma_k |\theta_k|, \qquad I_2(\theta | \tilde{\theta}) := I(\theta) - I_1(\theta | \tilde{\theta}).$$

Likewise for the empirical versions:

$$\hat{I}_1(\theta|\tilde{\theta}) := \sum_{k\,:\,\tilde{\theta}_k \neq 0} \hat{\sigma}_k |\theta_k|, \qquad \hat{I}_2(\theta|\tilde{\theta}) := \hat{I}(\theta) - \hat{I}_1(\theta|\tilde{\theta}).$$

## 2. Main results.

2.1. *Nonrandom normalization weights in the penalty.* Suppose that the $\sigma_k$ are known. Consider the estimator

$$(4) \qquad \hat{\theta}_n := \arg\min_{\theta \in \Theta} \{P_n \gamma_{f_\theta} + \lambda_n I(\theta)\}.$$

We now define the following six quantities:

(1) $\lambda_n := 2\bar{\lambda}_{n,0}$,
(2) $\mathcal{V}_\theta := H(4\lambda_n\sqrt{D_\theta})$,
(3) $\theta_n^* := \arg\min_{\theta \in \Theta}\{\mathcal{E}(f_\theta) + \mathcal{V}_\theta\}$,
(4) $2\epsilon_n^* := 3\mathcal{E}(f_{\theta_n^*}) + 2\mathcal{V}_{\theta_n^*}$,
(5) $\zeta_n^* := \epsilon_n^*/\bar{\lambda}_{n,0}$,
(6) $\theta(\epsilon_n^*) := \arg\min_{\theta \in \Theta, I(\theta - \theta_n^*) \leq 6\zeta_n^*}\{\mathcal{E}(f_\theta) - 4\lambda_n I_1(\theta - \theta_n^*|\theta_n^*)\}$.

These are defined "locally," that is, only for this subsection. This is because we have set some arbitrary values for certain constants involved. In Section 2.2, the six quantities will appear as well, but now with other constants. Moreover, in the Appendix, we define them explicitly as function of the constants.

CONDITION I. It holds that $\|f_{\theta_n^*} - \bar{f}\|_\infty \leq \eta$, where $\eta$ is given in Assumption B.

CONDITION II. It holds that $\|f_{\theta(\epsilon_n^*)} - \bar{f}\|_\infty \leq \eta$, where $\eta$ is given in Assumption B.

THEOREM 2.1. *Suppose Assumptions* L, A, B *and* C, *and Conditions* I *and* II *hold. Let* $\lambda_n$, $\theta_n^*$, $\epsilon_n^*$ *and* $\zeta_n^*$ *be given in* (1)–(6). *Assume* $\sigma_k$ *is known for all* $k$ *and let* $\hat{\theta}_n$ *be given in* (4). *Then we have with probability at least*

$$1 - 7\exp[-n\bar{a}_n^2 t^2],$$

*that*

$$(5) \qquad\qquad \mathcal{E}(f_{\hat{\theta}_n}) \leq 2\epsilon_n^*,$$

*and moreover*

$$(6) \qquad\qquad 2I(\hat{\theta}_n - \theta_n^*) \leq 7\zeta_n^*.$$

We now first discuss the quantities (1)–(6) and then the message of the theorem. To appreciate the meaning of the result, it may help to consider the "typical" case, with $G$ (given in Assumption B) quadratic [say $G(u) = u^2/2$], so that $H$ is also quadratic [say $H(v) = v^2/2$], and with $D(\mathcal{K})$ (given in Assumption C) being (up to constants) the cardinality $|\mathcal{K}|$ of the set index set $\mathcal{K}$ (see Section 3.1 for a discussion).

First, recall that $\lambda_n$ in (1) is the smoothing parameter. Note thus that $\lambda_n$ is chosen to be at least $(2\times)(4\times)\sqrt{2\log(2m)/n}$.

The function $\mathcal{V}_\theta$ in (2) depends only on the set of nonzero coefficients in $\theta$, and we refer to it as "estimation error" (in a generic sense). Note that in the "typical" case, $\mathcal{V}_\theta$ is up to constants equal to $\lambda_n^2 \dim_\theta$, where $\dim_\theta := |\{k : \theta_k \neq 0\}|$, that is, it is of order $\log m \times \dim_\theta/n$.

Because $\theta_n^*$ in (3) balances approximation error $\mathcal{E}(f_\theta)$ and "estimation error" $\mathcal{V}_\theta$, we refer to it as the "oracle" (although in different places, oracles may differ due to various choices of certain constants).

The terminology "estimation error" and "oracle" is inspired by results for the "typical" case, because for estimating a given number, $\dim_*$ say, of coefficients, the estimation error is typically of order $\dim_*/n$. The additional factor $\log m$ is the price for not knowing which coefficients are relevant.

Because we only show that $\mathcal{V}_\theta$ is an upper bound for estimation error, our terminology is not really justified in general. Moreover, our "oracle" is not allowed to try a different loss function with perhaps better margin behavior. To summarize, our terminology is mainly chosen for ease of reference.

The quantity $\epsilon_n^*$ in (4) will be called the "oracle rate" (for the excess risk). Moreover, $\zeta_n^*$ in (5) will be called the "oracle rate" for the $\ell_1$ norm. The latter can be compared with the result for Gaussian regression with orthogonal (and fixed) design (where $m \leq n$). Then the soft-thresholding estimator converges in $\ell_1$ norm with rate $\dim^* \sqrt{\log n/n}$ with $\dim^*$ up to constants the number of coefficients larger than $\sqrt{\log n/n}$ of the target. In the "typical" situation, $\zeta_n^*$ is of corresponding form.

The last quantity (6), the vector $\theta(\epsilon_n^*)$, balances excess risk with "being different" from $\theta_n^*$. Note that the balance is done locally, near values of $\theta$ where the $\ell_1$ distance between $\theta$ and $\theta_n^*$ is at most $6\zeta_n^*$.

Conditions I and II are technical conditions, which we need because Assumption B is only assumed to hold in a "neighborhood" of the target. Condition II may follow from a fast enough rate of convergence of the oracle.

The theorem states that the estimator with lasso penalty has excess risk at most twice $\epsilon_n^*$. It achieves this without knowledge about the estimation error $\mathcal{V}_\theta$ or, in particular, about the conjugate $H(\cdot)$ or the function $D(\cdot)$. In this sense, the estimator mimics the oracle.

The smoothing parameter $\lambda_n$, being larger than $\sqrt{2\log(2m)/n}$, plays its role in the estimation error $\mathcal{V}_{\theta_n^*}$. So from asymptotic point of view, we have the usual requirement that $m$ does not grow exponentially in $n$. In fact, in order for the estimation error to vanish, we require $D_{\theta_n^*} \log m/n \to 0$.

We observe that the constants involved are explicit, and reasonable for finite sample sizes. In the Appendix, we formulate the results with more general constants. For example, it is allowed there that the smoothing parameter $\lambda_n$ is not necessarily twice $\bar{\lambda}_{n,0}$, but may be arbitrary close to $\bar{\lambda}_{n,0}$, with consequences on the oracle rate. Our specific choice of the constants is merely to present the result in an transparent way.

The constant $a_n$ is, roughly speaking (when $K_m$ is not very large), the usual threshold that occurs when considering the Gaussian linear model with orthogonal design (and $m \le n$). The factor "4" appearing in front of it in our definition of $\bar{a}_n$, is due to the fact that we use a symmetrization inequality (which accounts for a factor "2") and a contraction inequality (which accounts for another factor "2"). We refer to Lemma A.2 for the details.

In our proof, we make use of Bousquet's inequality [Bousquet (2002)] for the amount of concentration of the empirical process around its mean. This inequality has the great achievement that the constants involved are economical. There is certainly room for improvement in the constants we, in our turn, derived from applying Bousquet's inequality. Our result should, therefore, only be seen as indication that the theory has something to say about finite sample sizes. In practice, cross validation for choosing the smoothing parameter seems to be the most reasonable way to go.

Finally, Theorem 2.1 shows that the $\ell_1$ difference between the estimated coefficients $\hat{\theta}_n$ and the oracle $\theta_n^*$ can be small. In this sense the lasso performs feature selection when the oracle rate is sufficiently fast. When $R(\bar{f})$ is not small, the convergence of $I(\hat{\theta}_n - \theta_n^*)$ is perhaps more important than fast rates for the excess risk.

2.2. *Random normalization weights in the penalty.* In this subsection, we estimate $\sigma_k^2 := Q\psi_k^2$ by $\hat{\sigma}_k^2 := Q_n\psi_k^2$, $k = 1, \ldots, m$. A complication is that the bound $K_m$ will generally also be unknown. The smoothing parameter $\lambda_n$ depends on $K_m$, as well as on $t$. We will assume a known but rough bound on $K_m$ (see Condition III' below). Then, we choose a known (large enough) value for $\bar{\lambda}_{n,0}$, which corresponds to an unknown value of $t$. This is elaborated upon in Theorem 2.2.

Recall the estimator

$$(7) \qquad \hat{\theta}_n = \underset{\theta \in \Theta}{\arg\min} \{ P_n \gamma_{f_\theta} + \lambda_n \hat{I}(\theta) \}.$$

In the new setup, we define the six quantities as:

$(1)'$  $\lambda_n := 3\bar{\lambda}_{n,0},$

$(2)'$  $\mathcal{V}_\theta := H(5\lambda_n\sqrt{D_\theta})$,
$(3)'$  $\theta_n^* := \arg\min_{\theta\in\Theta}\{\mathcal{E}(f_\theta)+\mathcal{V}_\theta\}$,
$(4)'$  $2\epsilon_n^* := 3\mathcal{E}(f_{\theta_n^*})+2\mathcal{V}_{\theta_n^*}$,
$(5)'$  $\zeta_n^* := \epsilon_n^*/\bar{\lambda}_{n,0}$.
$(6)'$  $\theta(\epsilon_n^*) := \arg\min_{\theta\in\Theta,\,I(\theta-\theta^*)\le 6\zeta_n^*}\{\mathcal{E}(f_\theta)-5\lambda_n I_1(\theta-\theta_n^*|\theta_n^*)\}$.

Note that $(1)'$–$(6)'$ are just (1)–(6) with slightly different constants.

CONDITION I$'$.   It holds that $\|f_{\theta_n^*}-\bar{f}\|_\infty \le \eta$, with $\eta$ given in Assumption B.

CONDITION II$'$.   It holds that $\|f_{\theta(\epsilon_n^*)}-\bar{f}\|_\infty \le \eta$, with $\eta$ given in Assumption B.

CONDITION III$'$.   We have $\sqrt{\frac{\log(2m)}{n}}K_m \le 0.13$.

The constant 0.13 in Condition III$'$ was again set quite arbitrary, in the sense that, provided some other constants are adjusted properly, it may be replaced by any other constant smaller than $(\sqrt{6}-\sqrt{2})/2$. The latter comes from our calculations using Bousquet's inequality.

THEOREM 2.2.   *Suppose Assumptions* L, A, B *and* C, *and Conditions* I$'$, II$'$ *and* III$'$ *are met. Let* $\lambda_n$, $\theta_n^*$, $\epsilon_n^*$ *and* $\zeta_n^*$ *be given in* $(1)'$–$(6)'$ *and let* $\hat{\theta}_n$ *be given in* (7). *Take*

$$\bar{\lambda}_{n,0} > 4\sqrt{\frac{\log(2m)}{n}} \times (1.6).$$

*Then with probability at least* $1-\alpha$, *we have that*

$$(8)\qquad\qquad\qquad \mathcal{E}(f_{\hat{\theta}_n}) \le 2\epsilon_n^*,$$

*and moreover*

$$(9)\qquad\qquad\qquad 2I(\hat{\theta}_n - \theta_n^*) \le 7\zeta_n^*.$$

*Here*

$$(10)\qquad\qquad \alpha = \exp[-na_n^2 s^2] + 7\exp[-n\bar{a}_n^2 t^2],$$

*with* $s>0$ *being defined by* $\frac{5}{9} = K_m\lambda_{n,0}(s)$, *and* $t>0$ *being defined by* $\bar{\lambda}_{n,0} = \bar{\lambda}_{n,0}(t)$.

The definition of $\lambda_{n,0}(s)$ $[\bar{\lambda}_{n,0}(t)]$ was given in (2) [(3)]. Thus, Theorem 2.2 gives qualitatively the same conclusion as Theorem 2.1.

To estimate the "confidence level" $\alpha$, we need an estimate of $K_m$. In Corollary A.3, we present an estimated upper bound for $\alpha$. An estimate of a lower bound can be found similarly.

## 3. Discussion of the assumptions and some examples.

3.1. *Discussion of Assumptions* L, A, B *and* C.  We assume in Assumption L that the loss function $\gamma$ is Lipschitz in $f$. This corresponds to the "robust case" with bounded influence function (e.g., Huber or quantile loss functions). The Lipschitz condition allows us to apply the contraction inequality (see Theorem A.3). It may well be that it depends on whether or not $\mathcal{F}$ is a class of functions uniformly bounded in sup norm. Furthermore, from the proofs one may conclude that we only need the Lipschitz condition locally near the target $\bar{f}$ [i.e., for those $\theta$ with $I(\theta - \theta_n^*) \leq 6\zeta_n^*$]. This will be exploited in Example 4 where the least squares loss is considered.

Assumption A is a technical but important condition. In fact, Condition III$'$ requires a bound proportional to $\sqrt{n/\log(2m)}$ for the constant $K_m$ of Assumption A. When for instance $\psi_1, \ldots, \psi_m$ actually form the co-variable $X$ itself [i.e., $X = (\psi_1(X), \ldots, \psi_m(X)) \in \mathbf{R}^m$], one possibly has that $K_m \leq K$ where $K$ does not grow with $m$. When the $\psi_k$ are truly feature mappings (e.g., wavelets), Assumption A may put a restriction on the number $m$ of base functions $\psi_k$.

Inspired by Tsybakov (2004), we call Assumption B the *margin* assumption. In most situations, $G$ is quadratic, that is, for some constant $C_0$, $G(u) = u^2/(2C_0)$, $u > 0$. This happens, for example, when $\mathbf{F}$ is the set, of *all* (measurable) functions, and in addition $\gamma_f(X, Y) = \gamma(f(X), Y)$ and $E(\gamma(z, Y)|X)$ is twice differentiable in a neighborhood of $z = \bar{f}(X)$, with second derivative at least $1/C_0$, $Q$-a.s.

Next, we address Assumption C. Define the vector $\psi = (\psi_1, \ldots, \psi_m)^T$. Suppose the $m \times m$ matrix

$$\Sigma := \int \psi \psi^T \, dQ,$$

has smallest eigenvalue $\beta^2 > 0$. Then one can easily verify that Assumption C holds with

$$(11) \qquad D(\mathcal{K}) = |\mathcal{K}|/\beta^2,$$

where $|\mathcal{K}|$ is the cardinality of the set $\mathcal{K}$. Thus, then $D$ indicates "dimension." Weighted versions are also relatively straightforward. Let $A(\mathcal{K})$ be the selection matrix:

$$\sum_{k \in \mathcal{K}} \alpha_k^2 = \alpha^T A(\mathcal{K}) \alpha, \qquad \alpha \in \mathbf{R}^m.$$

Then we can take

$$D(\mathcal{K}) = \sum_{k \in \mathcal{K}} w_k^{-2}/\beta^2(\mathcal{K}),$$

with $w_1, \ldots, w_m$ being a set of positive weights, and $1/\beta^2(\mathcal{K})$ the largest eigenvalue of the matrix $\Sigma^{-1/2} W A(\mathcal{K}) W \Sigma^{-1/2}$, $W := \operatorname{diag}(w_1, \ldots, w_m)$ [see also Tarigan and van de Geer (2006)].

A refinement of Assumption C, for example the *cumulative local coherence* assumption [see Bunea, Tsybakov and Wegkamp (2007a)] is needed to handle the case of overcomplete systems. We propose the refinement

ASSUMPTION C*. There exists nonnegative functions $\rho(\cdot)$ and $D(\cdot)$ on the subsets of the index set $\{1, \ldots, m\}$, such that for all $\mathcal{K} \subset \{1, \ldots, m\}$, and for all $\theta \in \Theta$ and $\tilde{\theta} \in \Theta$, we have

$$\sum_{k \in \mathcal{K}} \sigma_k |\theta_k - \tilde{\theta}_k| \leq \rho(\mathcal{K}) I(\theta - \tilde{\theta})/2 + \sqrt{D(\mathcal{K})} \|f_\theta - f_{\hat{\theta}}\|.$$

With Assumption C*, one can prove versions of Theorems 2.1 and 2.2 with different constants, provided that the "oracle" $\theta_n^*$ is restricted to having a value $\rho(\{K : \theta_{n,k}^* \neq 0\})$ strictly less than one. It can be shown that this is true under the cumulative local coherence assumption in Bunea, Tsybakov and Wegkamp (2007a), with $D(\mathcal{K})$ again proportional to $|\mathcal{K}|$. The reason why it works is because the additional term in Condition C* (as compared to Condition C) is killed by the penalty. Thus, the results can be extended to situations where $\Sigma$ is singular. However, due to lack of space we will not present a full account of this extension.

Assumption A is intertwined with Assumption C. As an illustration, suppose that one wants to include dummy variables for exclusive groups in the system $\{\psi_k\}$. Then Assumptions A and C together lead to requiring that the number of observations per group is not much smaller than $n/K_m^2$.

3.2. *Modifications.* In many cases, it is natural to leave a given subset of the coefficients not penalized, for example, those corresponding to the constant term and perhaps some linear terms, or to co-variables that are considered as definitely relevant. Suppose there are $p < m$ such coefficients, say the first $p$. The penalty is then modified to

$$\hat{I}(\theta) = \sum_{k=p+1}^m \hat{\sigma}_k |\theta_k|.$$

With this modification, the arguments used in the proof need only slight adjustments.

An important special case is where $\{\psi_k\}$ contains the constant function $\psi_1 \equiv 1$, and $\theta_1$ is not penalized. In that case, it is natural to modify Assumption A to

$$K_m := \max_{2 \leq k \leq m} \frac{\|\psi_k - \mu_k\|_\infty}{\sigma_k} < \infty,$$

where $\mu_k = Q\psi_k$ and where $\sigma_k^2$ is now defined as $\sigma_k^2 = Q\psi_k^2 - \mu_k^2$, $k = 2, \ldots, m$. Moreover, the penalty may be modified to $\hat{I}(\theta) = \sum_{k=2}^m \hat{\sigma}_k |\theta_k|$, where now $\hat{\sigma}_k^2 = Q_n\psi_k^2 - \hat{\mu}_k^2$ and $\hat{\mu}_k = Q_n\psi_k$. Thus, the means $\mu_k$ are now also estimated. However, this additional source of randomness is in a sense of smaller order. In conclusion, this modification does not bring in new theoretical complications, but does have a slight influence on the constants.

### 3.3. *Some examples.*

EXAMPLE 1 (Logistic regression). Consider the case $Y \in \{0, 1\}$ and the logistic loss function

$$\gamma_f(x, y) = \gamma(f(x), y) := [-f(x)y + \log(1 + e^{f(x)})]/2.$$

It is clear that this loss function is convex and Lipschitz. Let the target $\bar{f}$ be the log-odds ratio, that is $\bar{f} = \log(\frac{\pi}{1-\pi})$. where $\pi(x) := E(Y|X = x)$. It is easy to see that when for some $\varepsilon > 0$, it holds that $\varepsilon \leq \pi \leq 1 - \varepsilon$, $Q$-a.e., then Assumption B is met with $G(u) = u^2/(2C_0)$, $u > 0$, and with constant $C_0$ depending on $\eta$ and $\varepsilon$. The conjugate $H$ is then also quadratic, say $H(v) = C_1 v^2$, $v > 0$, where $C_1$ easily can be derived from $C_0$. Thus, for example, Theorem 2.2 has estimation error $\mathcal{V}_\theta = 25 C_1 \lambda_n^2 D_\theta$.

EXAMPLE 2 (Density estimation). Let $X$ have density $q_0 := dQ/d\nu$ with respect to a given $\sigma$-finite dominating measure $\nu$. We estimate this density by

$$\hat{q}_n = \exp[f_{\hat{\theta}_n} - b(f_{\hat{\theta}_n})],$$

where $\hat{\theta}_n$ is the lasso-penalized estimator with loss function

$$\gamma_f(X) := -f(X) + b(f), \qquad f \in \mathbf{F}.$$

Here $\mathbf{F} = \{f : \int e^f \, d\nu < \infty\}$ and $b(f) = \log \int e^f \, d\nu$ is the normalization constant. There is no response variable $Y$ in this case.

Clearly $\gamma_f$ satisfies Assumption L. Now, let $f_\theta = \sum_{k=1}^m \theta_k \psi_k$, with $\theta \in \Theta$ and $\Theta = \{\theta \in \mathbf{R}^m : f_\theta \in \mathbf{F}\}$. If we assume, for some $\varepsilon > 0$, that $\varepsilon \leq q_0 \leq 1/\varepsilon$, $\mu$-a.e, one easily verifies that Assumption B holds for $G(u) = u^2/(2C_0)$, $u > 0$, with $C_0$ depending on $\varepsilon$ and $\eta$. So we arrive at a similar estimation error as in the logistic regression example.

Now, a constant term will not be identifiable in this case, so we do not put the constant function in the system $\{\psi_k\}$. We may moreover take the functions $\psi_k$ centered in a convenient way, say $\int \psi_k \, d\nu = 0$, $\forall k$. In that situation a natural penalty with nonrandom weights could be $\sum_{k=1}^m w_k |\theta_k|$, where $w_k^2 = \int \psi_k^2 \, d\nu$, $k = 1, \ldots, m$. However, with this penalty, we are only able to prove a result along the lines of Theorem 2.1 or 2.2, when the smoothing parameter $\hat{\lambda}_n$ is chosen depending on an estimated lower bound for the density $q_0$. The penalty $\tilde{I}(\theta)$ with random weights $\hat{\sigma}_k = (Q_n \psi_k^2)^{1/2}$, $k = 1, \ldots, m$, allows a choice of $\lambda_n$ which does not depend on such an estimate.

Note that in this example, $\gamma(f(X), Y) = -f(X)$ is linear in $f$. Thus, in Lemma A.2, we do not need a contraction or symmetrization inequality, and may replace $\bar{a}_n$ $(\bar{\lambda}_{n,0})$ by $a_n$ $(\lambda_{n,0})$ throughout, that is, we may cancel a factor 4.

EXAMPLE 3 (Hinge loss). Let $Y \in \{\pm 1\}$. The hinge loss function, used in binary classification, is

$$\gamma_f(X, Y) := (1 - Yf(X))_+.$$

Clearly, hinge loss satisfies Assumption L. Let $\pi(x) := P(Y = 1|X = x)$ and **F** be *all* (measurable) functions. Then the target $\bar{f}$ is Bayes decision rule, that is $\bar{f} = \text{sign}(2\pi - 1)$. Tarigan and van de Geer (2006) show that instead of Condition B, it in fact holds, under reasonable conditions, that

$$\mathcal{E}(f_\theta) \geq G(\|f_\theta - \bar{f}\|_1^{1/2}),$$

with $G$ being again a strictly convex increasing function, but with $\| \cdot \|_1^{1/2}$ being the square root $L_1(Q)$ norm instead of the $L_2(Q)$ norm. The different norm however does not require essentially new theory. The bound $\|f - \bar{f}\| \leq \eta \|f - \bar{f}\|_1^{1/2}$ (which holds if $\|f - \bar{f}\|_\infty \leq \eta$) now ties Assumption B to Assumption C. See Tarigan and van de Geer (2006) for the details.

Alternatively, we may aim at a different target $\bar{f}$. Consider a minimizer $f_{\bar{\theta}}$ over $\mathcal{F} = \{f_\theta : \theta \in \Theta\}$, of expected hinge loss $E(1 - Yf(X))+$. Suppose that $\text{sign}(f_{\bar{\theta}})$ is Bayes decision rule $\text{sign}(2\pi - 1)$. Then $f_{\bar{\theta}}$ is a good target for classification. But for $\bar{f} = f_{\bar{\theta}}$, the margin behavior (Assumption B) is as yet not well understood.

EXAMPLE 4 (Quadratic loss). Suppose $Y = \bar{f}(X) + \varepsilon$, where $\varepsilon$ is $\mathcal{N}(0, 1)$ distributed and independent of $X$. Let $\gamma_f$ be quadratic loss

$$\gamma_f(X, Y) := \tfrac{1}{2}(Y - f(X))^2.$$

It is clear that in this case, Assumption B holds for all $f$, with $G(u) = u^2/2$, $u > 0$. However, the quadratic loss function is not Lipschitz on the whole real line. The Lipschitz property was used in order to apply the contraction inequality to the empirical process [see (15) of Lemma A.2 in Section 4.2]. To handle quadratic loss and Gaussian errors, we may apply exponential bounds for Gaussian random variables and a "local" Lipschitz property. Otherwise, one can use similar arguments as in the Lipschitz case. This leads to the following result for the case $\{\sigma_k\}$ known. (The case $\{\sigma_k\}$ is unknown can be treated similarly.)

THEOREM 3.1. *Suppose Assumptions* A *and* C *hold. Let* $\lambda_n$, $\theta_n^*$, $\epsilon_n^*$ *and* $\zeta_n^*$ *be given in* (1)–(6), *with* $H(v) = v^2/2$, $v > 0$, *but now with* $\bar{\lambda}_{n,0}$ *replaced by*

$$\tilde{\lambda}_{n,0} := \sqrt{\frac{14}{9}}\sqrt{\frac{2\log(2m)}{n} + 2t^2\bar{a}_n^2} + \bar{\lambda}_{n,0}.$$

*Assume moreover that* $\|f_{\theta_n^*} - \bar{f}\|_\infty \le \eta \le 1/2$, *that* $6\zeta_n^* K_m + 2\eta \le 1$, *and that* $\sqrt{\frac{\log(2m)}{n}}K_m \le 0.33$. *Let* $\sigma_k$ *be known for all* $k$ *and let* $\hat{\theta}_n$ *be given in* (4). *Then we have with probability at least* $1 - \alpha$, *that*

(12)
$$\mathcal{E}(f_{\hat{\theta}_n}) \le 2\epsilon_n^*,$$

*and moreover*

(13)
$$2I(\hat{\theta}_n - \theta_n^*) \le 7\zeta_n^*.$$

*Here*

$$\alpha = \exp[-na_n^2 s^2] + 7\exp[-n\bar{a}_n^2 t^2],$$

*with* $s > 0$ *a solution of* $\frac{9}{5} = K_m\lambda_{n,0}(s)$.

We conclude that when the oracle rate is small enough, the theory essentially goes through. For example, when $K_m = O(1)$ we require that the oracle has not more than $O(\sqrt{n/\log n})$ nonzero coefficients. This is in line with the results in Bunea, Tsybakov and Wegkamp (2007b). We also refer to Bunea, Tsybakov and Wegkamp (2007b) for possible extensions, including non-Gaussian errors and overcomplete systems.

## APPENDIX

The proofs of the results in this paper have elements that have become standard in the literature on penalized M-estimation. We refer to Massart (2000b) for a rather complete account of these. The $\ell_1$ penalty however has as special feature that it allows one to avoid estimating explicitly the estimation error. Moreover, a new element in the proof is the way we use the convexity of the penalized loss function to enter directly into local conditions.

The organization of the proofs is as follows. We start out in the next subsection with general results for empirical processes. This is applied in Section A.2 to obtain, for all $M > 0$, a bound for the empirical process uniformly over $\{\theta \in \Theta : I(\theta - \theta^*) \le M\}$. Here, $\theta^*$ is some fixed value, which can be chosen conveniently, according to the situation. We will choose $\theta^*$ to be the oracle $\theta_n^*$. Once we have this, we can start the iteration process to obtain a bound for $I(\hat{\theta}_n - \theta_n^*)$. Given this bound, we then proceed proving a bound for the excess risk $\mathcal{E}(\hat{f}_n)$. Section A.3 does this when the $\sigma_k$ are known, and Section A.4 goes through the adjustments when the $\sigma_k$ are estimated. Section A.5 considers quadratic loss.

Throughout Sections A.1–A.4, we require that Assumptions L, A, B and C hold.

**A.1. Preliminaries.** Our main tool is a result from Bousquet (2002), which is an improvement of the constants in Massart (2000a), the latter being again an improvement of the constants in Ledoux (1996). The result says that the supremum of any empirical process is concentrated near its mean. The amount of concentration depends only on the maximal sup norm and the maximal variance.

THEOREM A.1 (Concentration theorem [Bousquet (2002)]). *Let* $Z_1, \ldots, Z_n$ *be independent random variables with values in some space* $\mathcal{Z}$ *and let* $\Gamma$ *be a class of real-valued functions on* $\mathcal{Z}$, *satisfying for some positive constants* $\eta_n$ *and* $\tau_n$

$$\|\gamma\|_\infty \leq \eta_n \qquad \forall\ \gamma \in \Gamma$$

*and*

$$\frac{1}{n}\sum_{i=1}^{n} \operatorname{var}(\gamma(Z_i)) \leq \tau_n^2 \qquad \forall\ \gamma \in \Gamma.$$

*Define*

$$\mathbf{Z} := \sup_{\gamma \in \Gamma} \left| \frac{1}{n}\sum_{i=1}^{n}(\gamma(Z_i) - E\gamma(Z_i)) \right|.$$

*Then for* $z > 0$,

$$\mathbf{P}\left( \mathbf{Z} \geq \mathbf{E}\mathbf{Z} + z\sqrt{2(\tau_n^2 + 2\eta_n \mathbf{E}\mathbf{Z})} + \frac{2z^2\eta_n}{3} \right) \leq \exp[-nz^2].$$

Bousquet's inequality involves the expectation of the supremum of the empirical process. This expectation can be a complicated object, but one may derive bounds for it using symmetrization, and—in our Lipschitz case— contraction. To state these techniques, we need to introduce i.i.d. random variables $\epsilon_1, \ldots, \epsilon_n$, taking values $\pm 1$ each with probability $1/2$. Such a sequence is called a Rademacher sequence.

THEOREM A.2 (Symmetrization theorem [van der Vaart and Wellner (1996)]). *Let* $Z_1, \ldots, Z_n$ *be independent random variables with values in* $\mathcal{Z}$, *and let* $\varepsilon_1, \ldots, \varepsilon_n$ *be a Rademacher sequence independent of* $Z_1, \ldots, Z_n$. *Let* $\Gamma$ *be a class of real-valued functions on* $\mathcal{Z}$. *Then*

$$\mathbf{E}\left( \sup_{\gamma \in \Gamma} \left| \sum_{i=1}^{n}\{\gamma(Z_i) - E\gamma(Z_i)\} \right| \right) \leq 2\mathbf{E}\left( \sup_{\gamma \in \Gamma} \left| \sum_{i=1}^{n}\varepsilon_i\gamma(Z_i) \right| \right).$$

THEOREM A.3 (Contraction theorem [Ledoux and Talagrand (1991)]).
*Let $z_1, \ldots, z_n$ be nonrandom elements of some space $\mathcal{Z}$ and let $\mathcal{F}$ be a class of real-valued functions on $\mathcal{Z}$. Consider Lipschitz functions $\gamma_i \colon \mathbf{R} \to \mathbf{R}$, that is,*

$$|\gamma_i(s) - \gamma_i(\tilde{s})| \leq |s - \tilde{s}| \qquad \forall \ s, \tilde{s} \in \mathbf{R}.$$

*Let $\varepsilon_1, \ldots, \varepsilon_n$ be a Rademacher sequence. Then for any function $f^* \colon \mathcal{Z} \to \mathbf{R}$, we have*

$$\mathbf{E}\left(\sup_{f \in \mathcal{F}} \left| \sum_{i=1}^{n} \varepsilon_i \{\gamma_i(f(z_i)) - \gamma_i(f^*(z_i))\} \right| \right)$$

$$\leq 2\mathbf{E}\left(\sup_{f \in \mathcal{F}} \left| \sum_{i=1}^{n} \varepsilon_i (f(z_i) - f^*(z_i)) \right| \right).$$

Now, suppose $\Gamma$ is a finite set of functions. In that case, a bound for the expectation of the supremum of the empirical process over all $\gamma \in \Gamma$ can be derived from Bernstein's inequality.

LEMMA A.1. *Let $Z_1, \ldots, Z_n$ be independent $\mathcal{Z}$-valued random variables, and $\gamma_1, \ldots, \gamma_m$ be real-valued functions on $\mathcal{Z}$, satisfying for $k = 1, \ldots, m$,*

$$\mathbf{E}\gamma_k(Z_i) = 0, \forall \ i \qquad \|\gamma_k\|_\infty \leq \eta_n, \ \frac{1}{n}\sum_{i=1}^{n} \mathbf{E}\gamma_k^2(Z_i) \leq \tau_n^2.$$

*Then*

$$\mathbf{E}\left(\max_{1 \leq k \leq m} \left| \frac{1}{n}\sum_{i=1}^{n} \gamma_k(Z_i) \right| \right) \leq \sqrt{\frac{2\tau_n^2 \log(2m)}{n}} + \frac{\eta_n \log(2m)}{n}.$$

PROOF. Write $\bar{\gamma}_k := \frac{1}{n}\sum_{i=1}^{n} \gamma_k(Z_i)$, $k = 1, \ldots, m$. A classical intermediate step of the proof of Bernstein's inequality [see, e.g., van de Geer (2000)] tells us that for $n/\beta > \eta_n$, we have

$$\mathbf{E}\exp(\beta\bar{\gamma}_k) \leq \exp\left(\frac{\beta^2 \tau_n^2}{2(n - \beta\eta_n)}\right).$$

The same is true if we replace $\bar{\gamma}_k$ by $-\bar{\gamma}_k$. Hence,

$$\mathbf{E}\left(\max_{1 \leq k \leq m} |\bar{\gamma}_k| \right) \leq \frac{1}{\beta} \log\left(\mathbf{E}\exp\left(\beta \max_k \pm \bar{\gamma}_k\right)\right)$$

$$\leq \frac{\log(2m)}{\beta} + \frac{\beta\tau_n^2}{2(n - \beta\eta_n)}.$$

Now, take

$$\frac{n}{\beta} = \eta_n + \sqrt{\frac{n\tau_n^2}{2\log(2m)}}. \qquad \qquad \square$$

**A.2. First application to M-estimation with lasso penalty.** We now turn to our specific context. We let $\varepsilon_1, \ldots, \varepsilon_n$ be a Rademacher sequence, independent of the training set $(X_1, Y_1), \ldots, (X_n, Y_n)$. Moreover, we fix some $\theta^* \in \Theta$ and let for $M > 0$, $\mathcal{F}_M := \{f_\theta : \theta \in \Theta,\ I(\theta - \theta^*) \leq M\}$ and

$$(14) \qquad \mathbf{Z}(M) := \sup_{f \in \mathcal{F}_M} |(P_n - P)(\gamma_{f_\theta} - \gamma_{f_{\theta^*}})|,$$

where $\gamma_f(X, Y) = \gamma(f(X), Y) + b(f)$ now denotes the loss function.

LEMMA A.2. *We have*

$$\mathbf{EZ}(M) \leq 4M\mathbf{E}\left(\max_{1 \leq k \leq m} \left|\frac{1}{n}\sum_{i=1}^{n} \varepsilon_i \psi_k(X_i)/\sigma_k\right|\right).$$

PROOF. By the Symmetrization Theorem,

$$\mathbf{EZ}(M) \leq 2\mathbf{E}\left(\sup_{f \in \mathcal{F}_M} \left|\frac{1}{n}\sum_{i=1}^{n} \varepsilon_i\{\gamma(f_\theta(X_i), Y_i) - \gamma(f_{\theta^*}(X_i), Y_i)\}\right|\right).$$

Now, let $(\mathbf{X}, \mathbf{Y}) = \{(X_i, Y_i)\}_{i=1}^{n}$ denote the sample, and let $\mathbf{E}_{(\mathbf{X}, \mathbf{Y})}$ denote the conditional expectation given $(\mathbf{X}, \mathbf{Y})$. Then invoke the Lipschitz property of the loss functions, and apply the Contraction Theorem, with $z_i := X_i$ and $\gamma_i(\cdot) := \gamma(Y_i, \cdot)$, $i = 1, \ldots, n$, to find

$$(15) \quad \begin{aligned} &\mathbf{E}_{(\mathbf{X}, \mathbf{Y})}\left(\sup_{f \in \mathcal{F}_M} \left|\frac{1}{n}\sum_{i=1}^{n} \varepsilon_i\{\gamma(f_\theta(X_i), Y_i) - \gamma(f_{\theta^*}(X_i), Y_i)\}\right|\right) \\ &\leq 2\mathbf{E}_{(\mathbf{X}, \mathbf{Y})}\left(\sup_{f \in \mathcal{F}_M} \left|\frac{1}{n}\sum_{i=1}^{n} \varepsilon_i(f_\theta(X_i) - f_{\theta^*}(X_i))\right|\right). \end{aligned}$$

But clearly,

$$\begin{aligned} \left|\frac{1}{n}\sum_{i=1}^{n} \varepsilon_i(f_\theta(X_i) - f_{\theta^*}(X_i))\right| &\leq \sum_{k=1}^{m} \sigma_k |\theta_k - \theta_k^*| \max_{1 \leq k \leq m} \left|\frac{1}{n}\sum_{i=1}^{n} \varepsilon_i \psi_k(X_i)/\sigma_k\right| \\ &= I(\theta - \theta^*) \max_{1 \leq k \leq m} \left|\frac{1}{n}\sum_{i=1}^{n} \varepsilon_i \psi_k(X_i)/\sigma_k\right|. \end{aligned}$$

Since for $f_\theta \in \mathcal{F}_M$, we have $I(\theta - \theta^*) \leq M$, the result follows. □

Our next task is to bound the quantity

$$\mathbf{E}\left(\max_{1 \leq k \leq m} \frac{|1/n \sum_{i=1}^{n} \varepsilon_i \psi_k(X_i)|}{\sigma_k}\right).$$

Recall definition (1):

$$a_n = \left(\sqrt{\frac{2\log(2m)}{n}} + \frac{\log(2m)}{n}K_m\right).$$

Lemma A.3.   *We have*

$$\mathbf{E}\left(\max_{1\le k\le m}\left|\frac{(Q_n-Q)(\psi_k)}{\sigma_k}\right|\right)\le a_n,$$

*as well as*

$$\mathbf{E}\left(\max_{1\le k\le m}\frac{|1/n\sum_{i=1}^n\varepsilon_i\psi(X_i)|}{\sigma_k}\right)\le a_n.$$

Proof.   This follows from $\|\psi_k\|_\infty/\sigma_k\le K_m$ and $\mathrm{var}(\psi_k(X))/\sigma_k^2\le 1$. So we may apply Lemma A.1 with $\eta_n=K_m$ and $\tau_n^2=1$.   □

We now arrive at the result that fits our purposes.

Corollary A.1.   *For all $M>0$ and all $\theta\in\Theta$ with $I(\theta-\theta^*)\le M$, it holds that*

$$\|\gamma_{f_\theta}-\gamma_{f_{\theta^*}}\|_\infty\le MK_m$$

*and*

$$P(\gamma_{f_\theta}-\gamma_{f_{\theta^*}})^2\le M^2.$$

*Therefore, since by Lemmas A.2 and A.3, for all $M>0$,*

$$\frac{\mathbf{E}\mathbf{Z}(M)}{M}\le\bar a_n,\qquad \bar a_n=4a_n,$$

*we have, in view of Bousquet's Concentration theorem, for all $M>0$ and all $t>0$,*

$$\mathbf{P}\left(\mathbf{Z}(M)\ge\bar a_nM\left(1+t\sqrt{2(1+2\bar a_nK_m)}+\frac{2t^2\bar a_nK_m}{3}\right)\right)\le\exp[-n\bar a_n^2t^2].$$

**A.3. Proofs of the results in Section 2.1.**   In this subsection, we assume that $\sigma_k$ are known and we consider the estimator

$$\hat\theta_n:=\operatorname*{arg\,min}_{\theta\in\Theta}\left\{\frac{1}{n}\sum_{i=1}^n\gamma_{f_\theta}(X_i,Y_i)+\lambda_nI(\theta)\right\}.$$

Take $b>0$, $d>1$, and

$$d_b:=d\left(\frac{b+d}{(d-1)b}\vee 1\right).$$

Let

(A1)  $\lambda_n:=(1+b)\bar\lambda_{n,0}$
(A2)  $\mathcal{V}_\theta:=2\delta H(\frac{2\lambda_n\sqrt{\mathcal{D}_\theta}}{\delta})$, where $0<\delta<1$,

(A3)  $\theta_n^* := \arg\min_{\theta \in \Theta} \{\mathcal{E}(f_\theta) + \mathcal{V}_\theta\}$,

(A4)  $\epsilon_n^* := (1+\delta)\mathcal{E}(f_{\theta_n^*}) + \mathcal{V}_{\theta_n^*}$,

(A5)  $\zeta_n^* := \frac{\epsilon_n^*}{\lambda_{n,0}}$,

(A6)  $\theta(\epsilon_n^*) := \arg\min_{\theta \in \Theta, I(\theta - \theta_n^*) \le d_b \zeta_n^*/b} \{\delta\mathcal{E}(f_\theta) - 2\lambda_n I_1(\theta - \theta_n^*|\theta_n^*)\}$.

CONDITION I$(b,\delta)$.  It holds that $\|f_{\theta_n^*} - \bar{f}\|_\infty \le \eta$.

CONDITION II$(b,d,\delta)$.  It holds that $\|f_{\theta(\epsilon_n^*)} - \bar{f}\|_\infty \le \eta$.

In this subsection, we let $\lambda_n$, $\mathcal{V}_\theta$, $\theta_n^*$, $\epsilon_n^*$, $\zeta_n^*$ and $\theta(\epsilon_n^*)$ be defined in (A1)–(A6). Theorem 2.1 takes $b = 1$, $\delta = 1/2$ and $d = 2$.

We start out with proving a bound for the $\ell_1$ norm of the coefficients $\theta$, when restricted to the set where $\theta_{n,k}^* \neq 0$, in terms of the excess risk $\mathcal{E}(f_\theta)$.

LEMMA A.4.  *Suppose Condition* I$(b,\delta)$ *and Condition* II$(b,\delta,d)$ *are met. For all* $\theta \in \Theta$ *with* $I(\theta - \theta_n^*) \le d_b \zeta_n^*/b$, *it holds that*

$$2\lambda_n I_1(\theta - \theta_n^*|\theta_n^*) \le \delta\mathcal{E}(f_\theta) + \epsilon_n^* - \mathcal{E}(f_{\theta_n^*}).$$

PROOF.  We use the short-hand notation $I_1(\theta) = I_1(\theta|\theta_n^*)$, $\theta \in \Theta$. When $I(\theta - \theta_n^*) \le d_b \zeta_n^*/b$, one has

$$2\lambda_n I_1(\theta - \theta_n^*) = 2\lambda_n I_1(\theta - \theta_n^*) - \delta\mathcal{E}(f_\theta) + \delta\mathcal{E}(f_\theta)$$
$$\le 2\lambda_n I_1(\theta(\epsilon_n^*) - \theta_n^*) - \delta\mathcal{E}(f_{\theta(\epsilon^*)}) + \delta\mathcal{E}(f_\theta).$$

By Assumption C, combined with Condition II$(b,\delta,d)$,

$$2\lambda_n I_1(\theta(\epsilon_n^*) - \theta_n^*) \le 2\lambda_n \sqrt{D_{\theta_n^*}} \|f_{\theta(\epsilon_n^*)} - f_{\theta_n^*}\|.$$

By the triangle inequality,

$$2\lambda_n \sqrt{D_{\theta_n^*}} \|f_{\theta(\epsilon_n^*)} - f_{\theta_n^*}\| \le 2\lambda_n \sqrt{D_{\theta_n^*}} \|f_{\theta(\epsilon_n^*)} - \bar{f}\| + 2\lambda_n \sqrt{D_{\theta_n^*}} \|f_{\theta_n^*} - \bar{f}\|.$$

Since $\|f_{\theta(\epsilon_n^*)} - \bar{f}\| \le \eta$ as well as $\|f_{\theta_n^*} - \bar{f}\| \le \eta$, it follows from Condition I$(b,\delta)$ and Condition II$(b,\delta,d)$, combined with Assumption B, that

$$2\lambda_n \sqrt{D_{\theta_n^*}} \|f_{\theta(\epsilon_n^*)} - f_{\theta_n^*}\| \le \delta\mathcal{E}(f_{\theta(\epsilon_n^*)}) + \delta\mathcal{E}(f_{\theta_n^*}) + \mathcal{V}_{\theta_n^*}.$$

Hence, when $I(\theta - \theta_n^*) \le d_b \zeta_n^*/b$,

$$2\lambda_n I_1(\theta - \theta_n^*) \le \delta\mathcal{E}(f_\theta) + \delta\mathcal{E}(f_{\theta_n^*}) + \mathcal{V}_{\theta_n^*} = \delta\mathcal{E}(f_\theta) + \epsilon_n^* - \mathcal{E}(f_{\theta_n^*}). \qquad \square$$

We now show that for any $\tilde{\theta}$ for which the penalized empirical risk is not larger than that of the oracle $\theta_n^*$, the $\ell_1$ difference between $\tilde{\theta}$ and $\theta_n^*$ has to be "small." Lemma A.5 represents one iteration, which shows that on the set $I(\tilde{\theta} - \theta_n^*) \leq d_0 \zeta_n^*/b$, in fact, except on a subset with small probability, $I(\tilde{\theta} - \theta_n^*)$ is strictly smaller than $d_0 \zeta_n^*/b$.

Let for $M > 0$, $\mathbf{Z}(M)$ be the random variable defined in (14), with $\theta^* = \theta_n^*$.

LEMMA A.5.   *Suppose Condition* I$(b, \delta)$ *and Condition* II$(b, \delta, d)$ *are met. Consider any* (*random*) $\tilde{\theta} \in \Theta$ *with* $R_n(f_{\tilde{\theta}}) + \lambda_n I(\tilde{\theta}) \leq R_n(f_{\theta_n^*}) + \lambda_n I(\theta_n^*)$. *Let* $1 < d_0 \leq d_b$. *Then*

$$\mathbf{P}\left( I(\tilde{\theta} - \theta_n^*) \leq d_0 \frac{\zeta_n^*}{b} \right)$$

$$\leq \mathbf{P}\left( I(\tilde{\theta} - \theta_n^*) \leq \left( \frac{d_0 + b}{1 + b} \right) \frac{\zeta_n^*}{b} \right) + \exp[-n \bar{a}_n^2 t^2].$$

PROOF.   We let $\tilde{\mathcal{E}} := \mathcal{E}(f_{\tilde{\theta}})$ and $\mathcal{E}^* := \mathcal{E}(f_{\theta_n^*})$. We also use the short hand notation $I_1(\theta) = I_1(\theta | \theta_n^*)$ and $I_2(\theta) = I_2(\theta | \theta_n^*)$. Since $R_n(f_{\tilde{\theta}}) + \lambda_n I(\tilde{\theta}) \leq R_n(f_{\theta_n^*}) + \lambda_n I(\theta_n^*)$, we know that when $I(\tilde{\theta} - \theta_n^*) \leq d_0 \zeta_n^*/b$, that

$$\tilde{\mathcal{E}} + \lambda_n I(\tilde{\theta}) \leq \mathbf{Z}(d_0 \zeta_n^*/b) + \mathcal{E}^* + \lambda_n I(\theta_n^*).$$

With probability at least $1 - \exp[-n \bar{a}_n^2 t^2]$, the random variable $\mathbf{Z}(d_0 \zeta_n^*/b)$ is bounded by $\bar{\lambda}_{n,0} d_0 \zeta_n^*/b$. But then we have

$$\tilde{\mathcal{E}} + \lambda_n I(\tilde{\theta}) \leq \bar{\lambda}_{n,0} \frac{d_0 \zeta_n^*}{b} + \mathcal{E}^* + \lambda_n I(\theta_n^*).$$

Invoking $\lambda_n = (1 + b) \bar{\lambda}_{n,0}$, $I(\tilde{\theta}) = I_1(\tilde{\theta}) + I_2(\tilde{\theta})$ and $I(\theta_n^*) = I_1(\theta_n^*)$, we find on $\{ I(\tilde{\theta} - \theta_n^*) \leq d_0 \zeta_n^*/b \} \cup \{ \mathbf{Z}(d_0 \zeta_n^*/b) \leq \bar{\lambda}_{n,0} d_0 \zeta_n^*/b \}$, that

$$\tilde{\mathcal{E}} + (1 + b) \bar{\lambda}_{n,0} I_2(\tilde{\theta})$$

$$\leq \bar{\lambda}_{n,0} \frac{d_0 \zeta_n^*}{b} + \mathcal{E}^* + (1 + b) \bar{\lambda}_{n,0} I_1(\theta_n^*) - (1 + b) \bar{\lambda}_{n,0} I_1(\tilde{\theta})$$

$$\leq \bar{\lambda}_{n,0} \frac{d_0 \zeta_n^*}{b} + \mathcal{E}^* + (1 + b) \bar{\lambda}_{n,0} I_1(\tilde{\theta} - \theta_n^*).$$

But $I_2(\tilde{\theta}) = I_2(\tilde{\theta} - \theta_n^*)$. So if we add another $(1 + b) \bar{\lambda}_{n,0} I_1(\tilde{\theta} - \theta_n^*)$ to both left- and right-hand side of the last inequality, we obtain

$$\tilde{\mathcal{E}} + (1 + b) \bar{\lambda}_{n,0} I(\tilde{\theta} - \theta_n^*) \leq \bar{\lambda}_{n,0} \frac{d_0 \zeta_n^*}{b} + 2(1 + b) \bar{\lambda}_{n,0} I_1(\tilde{\theta} - \theta_n^*).$$

Since $d_0 \leq d_b$, we know from Lemma A.4 that this implies

$$\tilde{\mathcal{E}} + (1 + b) \bar{\lambda}_{n,0} I(\tilde{\theta} - \theta_n^*) \leq \bar{\lambda}_{n,0} \frac{d_0 \zeta_n^*}{b} + \delta \tilde{\mathcal{E}} + \epsilon_n^*$$

$$= (d_0 + b) \bar{\lambda}_{n,0} \frac{\zeta_n^*}{b} + \delta \tilde{\mathcal{E}},$$

as $\epsilon_n^* = \bar{\lambda}_{n,0}\zeta_n^*$. But then

$$(1-\delta)\tilde{\mathcal{E}} + (1+b)\bar{\lambda}_{n,0}I(\tilde{\theta}-\theta_n^*) \le (d_0+b)\bar{\lambda}_{n,0}\frac{\zeta_n^*}{b}.$$

Because $0 < \delta < 1$, the result follows. $\square$

One may repeat the argument of the previous lemma $N$ times, to get the next corollary.

COROLLARY A.2. *Suppose Condition* I$(b,\delta)$ *and Condition* II$(b,\delta,d)$ *are met. Let $d_0 \le d_b$. For any (random) $\tilde{\theta} \in \Theta$ with $R_n(f_{\tilde{\theta}}) + \lambda_n I(\tilde{\theta}) \le R_n(f_{\theta_n^*}) + \lambda_n I(\theta_n^*)$,*

$$\mathbf{P}\left(I(\tilde{\theta} - \theta_n^*) \le d_0\frac{\zeta_n^*}{b}\right)$$

$$\le \mathbf{P}\left(I(\tilde{\theta} - \theta_n^*) \le (1 + (d_0-1)(1+b)^{-N})\frac{\zeta_n^*}{b}\right) + N\exp[-n\bar{a}_n^2 t^2].$$

The next lemma considers a convex combination of $\hat{\theta}_n$ and $\theta_n^*$, such that the $\ell_1$ distance between this convex combination and $\theta_n^*$ is small enough.

LEMMA A.6. *Suppose Condition* I$(b,\delta)$ *and Condition* II$(b,\delta,d)$ *are met. Define*

$$\tilde{\theta}_s = s\hat{\theta}_n + (1-s)\theta_n^*,$$

*with*

$$s = \frac{d\zeta_n^*}{d\zeta_n^* + bI(\hat{\theta}_n - \theta_n^*)}.$$

*Then, for any integer $N$, with probability at least $1 - N\exp[-n\bar{a}_n^2 t^2]$ we have*

$$I(\tilde{\theta}_s - \theta_n^*) \le (1 + (d-1)(1+b)^{-N})\frac{\zeta_n^*}{b}.$$

PROOF. The loss function and penalty are convex, so

$$R_n(f_{\hat{\theta}_s}) + \lambda_n I(\tilde{\theta}_s) \le sR_n(f_{\hat{\theta}_n}) + s\lambda_n I(\hat{\theta}_n) + (1-s)R_n(f_{\theta_n^*}) + (1-s)\lambda_n I(\theta_n^*)$$

$$\le R_n(f_{\theta_n^*}) + \lambda_n I(\theta_n^*).$$

Moreover,

$$I(\tilde{\theta}_s - \theta_n^*) = sI(\hat{\theta}_n - \theta_n^*) = \frac{d\zeta_n^* I(\hat{\theta}_n - \theta_n^*)}{d\zeta_n^* + bI(\hat{\theta}_n - \theta_n^*)} \le d\frac{\zeta_n^*}{b}.$$

By the definition of $d_b$, we know $d \le d_b$. Now apply the previous lemma with $\tilde{\tilde{\theta}} = \tilde{\theta}_s$ and $d_0 = d$.  □

One may now deduce that not only the convex combination, but also $\hat{\theta}_n$ itself is close to $\theta_n^*$, in $\ell_1$ distance.

LEMMA A.7.  *Suppose Condition* I$(b, \delta)$ *and Condition* II$(b, \delta, d)$ *are met. Let* $N_1 \in \mathbf{N}$ *and* $N_2 \in \mathbf{N} \cup \{0\}$. *Define* $\delta_1 = (1+b)^{-N_1}$ $(N_1 \ge 1)$, *and* $\delta_2 = (1+b)^{-N_2}$. *With probability at least* $1 - (N_1 + N_2) \exp[-n\bar{a}_n^2 t^2]$, *we have*

$$I(\hat{\theta}_n - \theta_n^*) \le d(\delta_1, \delta_2) \frac{\zeta_n^*}{b},$$

*with*

$$d(\delta_1, \delta_2) = 1 + \left( \frac{1 + (d^2 - 1)\delta_1}{(d-1)(1-\delta_1)} \right) \delta_2.$$

PROOF.  We know from Lemma A.6 that with probability at least $1 - N_1 \exp[-n\bar{a}_n^2 t^2]$,

$$I(\tilde{\theta}_s - \theta_n^*) \le (1 + (d-1)\delta_1) \frac{\zeta_n^*}{b}.$$

But then

$$I(\hat{\theta}_n - \theta_n^*) \le \frac{d(1 + (d-1)\delta_1)}{(d-1)(1-\delta_1)} \frac{\zeta_n^*}{b} := \bar{d}_0 \frac{\zeta_n^*}{b}.$$

We have $\bar{d}_0 \le d_b$, since $0 < \delta_1 \le 1/(1+b)$, and by the definition of $d_b$. Therefore, we may Corollary A.2 to $\hat{\theta}_n$, with $d_0$ replaced by $\bar{d}_0$. We find that with probability at least $1 - (N_1 + N_2) \exp[-n\bar{a}_n^2 t^2]$,

$$I(\hat{\theta}_n - \theta_n^*) \le (1 + (d_0 - 1)\delta_2) \frac{\zeta_n^*}{b}$$

$$= \left( 1 + \left( \frac{1 + (d^2 - 1)\delta_1}{(d-1)(1-\delta_1)} \right) \delta_2 \right) \frac{\zeta_n^*}{b}. \qquad □$$

Recall the notation

$$(16) \qquad d(\delta_1, \delta_2) := 1 + \left( \frac{1 + (d^2 - 1)\delta_1}{(1-d)(1-\delta_1)} \right) \delta_2.$$

Write

$$(17) \qquad \Delta(b, \delta, \delta_1, \delta_2) := d(\delta_1, \delta_2) \frac{1 - \delta^2}{\delta b} \vee 1.$$

THEOREM A.4. *Suppose Condition* I$(b, \delta)$ *and Condition* II$(b, \delta, d)$ *are met. Let* $\delta_1$ *and* $\delta_2$ *be as in Lemma* A.7. *We have with probability at least*

$$1 - \left( \log_{1+b} \frac{(1+b)^2 \Delta(b, \delta, \delta_1, \delta_2)}{\delta_1 \delta_2} \right) \exp[-n\bar{a}_n^2 t^2],$$

*that*

$$\mathcal{E}(f_{\hat{\theta}_n}) \leq \frac{\epsilon_n^*}{1 - \delta},$$

*and moreover*

$$I(\hat{\theta}_n - \theta_n^*) \leq d(\delta_1, \delta_2) \frac{\zeta_n^*}{b}.$$

PROOF. Define $\hat{\mathcal{E}} := \mathcal{E}(\hat{f}_{\hat{\theta}_n})$ and $\mathcal{E}^* := \mathcal{E}(f_{\theta_n^*})$. We also again use the short hand notation $I_1(\theta) = I_1(\theta | \theta_n^*)$ and $I_2(\theta) = I_2(\theta | \theta_n^*)$. Set

$$c := \frac{\delta b}{1 - \delta^2}.$$

We consider the cases (a) $c < d(\delta_1, \delta_2)$ and (b) $c \geq d(\delta_1, \delta_2)$.

(a) Suppose first that $c < d(\delta_1, \delta_2)$. Let $J$ be an integer satisfying $(1 + b)^{J-1} c \leq d(\delta_1, \delta_2)$ and $(1+b)^J c > d(\delta_1, \delta_2)$. We consider the cases (a1) $c\zeta_n^*/b < I(\hat{\theta}_n - \theta_n^*) \leq d(\delta_1, \delta_2)\zeta_n^*/b$ and (a2) $I(\hat{\theta}_n - \theta_n^*) \leq c\zeta_n^*/b$.

(a1) If $c\zeta_n^*/b < I(\hat{\theta}_n - \theta_n^*) \leq d(\delta_1, \delta_2)\zeta_n^*/b$, then

$$(1+b)^{j-1} c\zeta_n^*/b < I(\hat{\theta}_n - \theta_n^*) \leq (1+b)^j c\zeta_n^*/b,$$

for some $j \in \{1, \dots, J\}$. Except on set with probability at most $\exp[-n\bar{a}_n^2 t^2]$, we thus have that

$$\hat{\mathcal{E}} + (1+b)\bar{\lambda}_{n,0} I(\hat{\theta}_n) \leq (1+b)\bar{\lambda}_{n,0} I(\hat{\theta}_n - \theta_n^*) + \mathcal{E}^* + (1+b)\bar{\lambda}_{n,0} I(\theta_n^*).$$

So then, by similar arguments as in the proof of Lemma A.5, we find

$$\hat{\mathcal{E}} \leq 2(1+b)\bar{\lambda}_{n,0} I_1(\hat{\theta}_n - \theta_n^*) + \mathcal{E}^*.$$

Since $d(\delta_1, \delta_2) \leq d_b$, we obtain $\hat{\mathcal{E}} \leq \epsilon_n^* + \delta\hat{\mathcal{E}}$, so then $\hat{\mathcal{E}} \leq \frac{\epsilon_n^*}{1-\delta}$.

(a2) If $I(\hat{\theta}_n - \theta_n^*) \leq c\zeta_n^*/b$, we find, except on a set with probability at most $\exp[-n\bar{a}_n^2 t^2]$, that

$$(18) \qquad \hat{\mathcal{E}} + (1+b)\bar{\lambda}_{n,0} I(\hat{\theta}_n) \leq \left( \frac{\delta}{1-\delta^2} \right) \bar{\lambda}_{n,0}\zeta_n^* + \mathcal{E}^* + (1+b)I(\theta_n^*),$$

which gives

$$\hat{\mathcal{E}} \leq \left( \frac{\delta}{1-\delta^2} \right) \bar{\lambda}_{n,0}\zeta_n^* + \mathcal{E}^* + (1+b)\bar{\lambda}_{n,0} I_1(\hat{\theta}_n - \theta_n^*)$$

$$\leq \left(\frac{\delta}{1-\delta^2}\right)\bar{\lambda}_{n,0}\zeta_n^* + \mathcal{E}^* + \frac{\delta}{2}\mathcal{E}^* + \frac{\mathcal{V}_{\theta_n^*}}{2} + \frac{\delta}{2}\hat{\mathcal{E}}$$

$$= \left(\frac{\delta}{1-\delta^2} + \frac{1}{2}\right)\epsilon_n^* + \frac{\mathcal{E}^*}{2} + \frac{\delta}{2}\hat{\mathcal{E}}$$

$$\leq \left(\frac{\delta}{1-\delta^2} + \frac{1}{2} + \frac{1}{2(1+\delta)}\right)\epsilon_n^*.$$

But this yields

$$\hat{\mathcal{E}} \leq \frac{2}{2-\delta}\left(\frac{\delta}{1-\delta^2} + \frac{1}{2} + \frac{1}{2(1+\delta)}\right)\epsilon_n^* = \frac{1}{1-\delta}\epsilon_n^*.$$

Furthermore, by Lemma A.7, we have with probability at least $1-(N_1+N_2)\exp[-n\bar{a}_n t^2]$, that

$$I(\hat{\theta}_n - \theta_n^*) \leq \frac{d(\delta_1, \delta_2)}{b}\zeta_n^*.$$

The result now follows from

$$J + 1 \leq \log_{1+b}\left(\frac{(1+b)^2 d(\delta_1, \delta_2)}{c}\right)$$

and

$$N_1 = \log_{1+b}\left(\frac{1}{\delta_1}\right), \qquad N_2 = \log_{1+b}\left(\frac{1}{\delta_2}\right).$$

(b) Finally, consider the case $c \geq d(\delta_1, \delta_2)$. Then on the set where $I(\hat{\theta}_n - \theta_n^*) \leq d(\delta_1, \delta_2)\zeta_n^*/b$, we again have that, except on a subset with probability at most $\exp[-n\bar{a}_n t^2]$,

$$\hat{\mathcal{E}} + (1+b)\bar{\lambda}_{n,0}I(\hat{\theta}_n) \leq d(\delta_1, \delta_2)\frac{\zeta_n^*}{b} + \mathcal{E}^* + (1+b)I(\theta_n^*)$$

$$\leq \left(\frac{\delta}{1-\delta^2}\right)\bar{\lambda}_{n,0}\zeta_n^* + \mathcal{E}^* + (1+b)I(\theta_n^*),$$

as

$$d(\delta_1, \delta_2) \leq c = \frac{\delta b}{1-\delta^2}.$$

So we arrive at the same inequality as (18) and we may proceed as there. Note finally that also in this case

$$(N_1 + N_2 + 1) \leq (N_1 + N_2 + 2)$$

$$= \log_{1+b}\frac{(1+b)^2}{\delta_1\delta^2}$$

$$= \log_{1+b}\frac{(1+b)^2\Delta(b, \delta, \delta_1, \delta_2)}{\delta_1\delta_2}.$$

□

PROOF OF THEOREM 2.1.   Theorem 2.1 is a special case of Theorem A.4, with $b = 1$, $\delta = 1/2$, $d = 2$ and $\delta_1 = \delta_2 = 1/2$.   □

**A.4. Proof of the results in Section 2.2.**   Let

$$\Omega = \{\sigma_k/c_1 \leq \hat{\sigma}_k \leq c_2\sigma_k \; \forall \; k\},$$

where $c_1 > 1$ and $c_2 = \sqrt{2c_1^2 - 1}/c_1$. We show that for some $s$ depending on $c_1$, the set $\Omega$ has probability at least $1 - \exp[-na_n^2 s^2]$.

LEMMA A.8.   *We have*

$$\mathbf{E}\left(\max_{1 \leq k \leq m}\left|\frac{\hat{\sigma}_k^2}{\sigma_k^2} - 1\right|\right) \leq a_n K_m.$$

PROOF.   Apply Lemma A.1 with $\eta_n = K_m^2$ and $\tau_n^2 = K_m^2$.   □

Recall now that for $s > 0$,

$$\lambda_{n,0}(s) := a_n\left(1 + s\sqrt{2(1 + 2a_n K_m)} + \frac{2s^2 a_n K_m}{3}\right).$$

LEMMA A.9.   *Let* $c_0 > 1$. *Suppose that*

$$2\sqrt{\frac{\log(2m)}{n}}K_m \leq \frac{\sqrt{6c_0^2 - 4} - \sqrt{2c_0^2}}{c_0}.$$

*Let* $c_1 > c_0$. *Then there exists a solution* $s > 0$ *of*

$$1 - \frac{1}{c_1^2} = K_m\lambda_{n,0}(s). \tag{19}$$

*Moreover, with this value of* $s$,

$$\mathbf{P}(\Omega) \geq 1 - \exp[-na_n^2 s^2].$$

PROOF.   We first note that

$$2\sqrt{\frac{\log(2m)}{n}}K_m \leq \frac{\sqrt{6c_0^2 - 4} - \sqrt{2c_0^2}}{c_0},$$

implies

$$a_n K_m \leq 1 - \frac{1}{c_0^2}.$$

Therefore, there exists a solution $s > 0$ satisfying (19). By Bousquet's inequality and Lemma A.8,

$$\mathbf{P}\left(\max_k \left|\frac{\hat{\sigma}_k^2}{\sigma_k^2} - 1\right| \geq K_m \lambda_{n,0}(s)\right) \leq \exp[-na_n^2 s^2].$$

In other words,

$$\mathbf{P}\left(\max_k \left|\frac{\hat{\sigma}_k^2}{\sigma_k^2} - 1\right| \geq 1 - \frac{1}{c_1^2}\right) \leq \exp[-na_n^2 s^2].$$

But the inequality

$$\left|\frac{\hat{\sigma}_k^2}{\sigma_k^2} - 1\right| \geq 1 - \frac{1}{c_1^2} \qquad \forall\, k$$

is equivalent to

$$\sigma_k/c_1 \leq \hat{\sigma} \leq c_2 \sigma_k \qquad \forall\, k. \qquad \square$$

Recall the estimator

$$\hat{\theta}_n = \arg\min_{\theta \in \Theta}\{R_n(f_\theta) + \lambda_n \hat{I}(\theta)\}.$$

In this case, for $1 + b > c_1 > 1$, $c_2 = \sqrt{2c_1^2 - 1}/c_1$, and $d > 1$, and for

$$d_b := d\left(\frac{d+b}{(d-1)b} \vee 1\right),$$

given as before, we define:

(A1)′  $\lambda_n := c_1(1+b)\bar{\lambda}_{n,0},$
(A2)′  $\mathcal{V}_\theta := 2\delta H\left(\frac{2c_2\lambda_n\sqrt{D_\theta}}{\delta}\right)$, where $0 < \delta < 1$,
(A3)′  $\theta_n^* := \arg\min_{\theta \in \Theta}\{\mathcal{E}(f_\theta) + \mathcal{V}_\theta\},$
(A4)′  $\epsilon_n^* := (1+\delta)\mathcal{E}(f_{\theta_n^*}) + \mathcal{V}_{\theta_n^*},$
(A5)′  $\zeta_n^* := \epsilon_n^*/\bar{\lambda}_{n,0},$
(A6)′  $\theta(\epsilon_n^*) := \arg\min_{\theta \in \Theta, I(\theta-\theta^*) \leq d_b \zeta_n^*/b}\{\delta\mathcal{E}(f_\theta) - 2c_2\lambda_n I_1(\theta - \theta_n^*|\theta_n^*)\}.$

Note that on $\Omega$,

$$I(\theta)/c_1 \leq \hat{I}(\theta) \leq c_2 I(\theta) \qquad \forall\, \theta$$

and

$$I_k(\theta|\theta_n^*)/c_1 \leq \hat{I}_k(\theta|\theta_n^*) \leq c_2 I_k(\theta|\theta_n^*), \qquad k = 1, 2, \ \forall\, \theta.$$

Hence, we have upper and lower bounds for the estimated $\ell_1$ norms. This means we may proceed as in the previous subsection.

CONDITION I($\delta, b, c_1, c_2$).  It holds that

$$\|f_{\theta_n^*} - \bar{f}\|_\infty \leq \eta.$$

CONDITION II($\delta, b, d, c_1, c_2$).   It holds that $\|f_{\theta(\epsilon_n^*)} - \bar{f}\|_\infty \leq \eta.$

CONDITION III($c_0$).  For some known constant $c_0 > 1$, it holds that $2\sqrt{\log(2m)/n} K_m \leq (\sqrt{6c_0^2 - 4} - \sqrt{2c_0^2})/c_0.$

LEMMA A.10.  *Suppose Condition* I($\delta, b, c_1, c_2$) *and Condition* II($\delta, b, d,$ $c_1, c_2$). *Let* $c_1 b/(1 + b - c_1) < d_0 \leq \bar{d}_b$. *For any (random)* $\tilde{\theta} \in \Theta$ *with* $R_n(f_{\tilde{\theta}}) + \lambda_n I(\tilde{\theta}) \leq R_n(f_{\theta_n^*}) + \lambda_n I(\theta_n^*)$, *we have*

$$\mathbf{P}\left(I(\tilde{\theta} - \theta_n^*) \leq d_0 \frac{\zeta_n^*}{b}\right)$$

$$\leq \mathbf{P}\left(I(\tilde{\theta} - \theta_n^*) \leq \left(\frac{c_1(d_0 + b)}{1 + b}\right)\frac{\zeta_n^*}{b}\right)$$

$$+ 1 - \mathbf{P}(\Omega) + \exp[n\bar{a}_n^2 t^2].$$

PROOF.  This is essentially repeating the argument of Lemma A.5. Let $\tilde{\mathcal{E}} := \mathcal{E}(f_{\tilde{\theta}})$ and $\mathcal{E}^* := \mathcal{E}(f_{\theta_n^*})$, and use the short hand notation $I_1(\theta) = I_1(\theta|\theta_n^*)$ and $I_2(\theta) = I_2(\theta|\theta_n^*)$. Likewise, $\hat{I}_1(\theta) = \hat{I}_1(\theta|\theta_n^*)$ and $\hat{I}_2(\theta) = \hat{I}_2(\theta|\theta_n^*)$.

On the set $\{I(\tilde{\theta} - \theta_n^*) \leq d_0\zeta_n^*/b\}$, we have, except on a subset with probability at most $\exp[-n\bar{a}_n^2 t^2]$, that

$$\tilde{\mathcal{E}} + \lambda_n \hat{I}(\tilde{\theta}) \leq \bar{\lambda}_{n,0} \frac{d_0 \zeta_n^*}{b} + \mathcal{E}^* + \lambda_n \hat{I}(\theta_n^*).$$

In the remainder of the proof, we will, therefore, only need to consider the set

$$\Omega \cap \{I(\tilde{\theta} - \theta_n^*) \leq d_0\zeta_n^*/b\} \cap \left\{\tilde{\mathcal{E}} + \lambda_n \hat{I}(\tilde{\theta}) \leq \bar{\lambda}_{n,0} \frac{d_0 \zeta_n^*}{b} + \mathcal{E}^* + \lambda_n \hat{I}(\theta_n^*)\right\}.$$

Invoking $\lambda_n = c_1(1 + b)\bar{\lambda}_{n,0}$, $\hat{I}(\tilde{\theta}) = \hat{I}_1(\tilde{\theta}) + \hat{I}_2(\tilde{\theta})$ and $\hat{I}(\theta_n^*) = \hat{I}_1(\theta_n^*)$, we find

$$\tilde{\mathcal{E}} + c_1(1 + b)\bar{\lambda}_{n,0}\hat{I}_2(\tilde{\theta})$$

$$\leq \bar{\lambda}_{n,0} \frac{d_0 \zeta_n^*}{b} + \mathcal{E}^* + c_1(1 + b)\bar{\lambda}_{n,0}\hat{I}_1(\theta_n^*) - c_1(1 + b)\bar{\lambda}_{n,0}\hat{I}_1(\tilde{\theta})$$

$$\leq \bar{\lambda}_{n,0} \frac{d_0 \zeta_n^*}{b} + \mathcal{E}^* + c_1(1 + b)\bar{\lambda}_{n,0}\hat{I}_1(\tilde{\theta} - \theta_n^*).$$

We add another $c_1(1 + b)\bar{\lambda}_{n,0}\hat{I}_1(\tilde{\theta} - \theta_n^*)$ to both left- and right-hand side of the last inequality, and then apply $\hat{I}_1(\tilde{\theta} - \theta_n^*) \leq c_2 I_1(\tilde{\theta} - \theta_n^*)$, to obtain

$$\tilde{\mathcal{E}} + c_1(1 + b)\bar{\lambda}_{n,0}\hat{I}(\tilde{\theta} - \theta_n^*)$$

$$\leq \bar{\lambda}_{n,0}\frac{d_0\zeta_n^*}{b} + \mathcal{E}^* + 2c_1(1+b)\bar{\lambda}_{n,0}\hat{I}_1(\tilde{\theta} - \theta_n^*)$$

$$\leq \bar{\lambda}_{n,0}\frac{d_0\zeta_n^*}{b} + \mathcal{E}^* + 2c_1(1+b)c_2\bar{\lambda}_{n,0}I_1(\tilde{\theta} - \theta_n^*).$$

Now $I(\tilde{\theta} - \theta_n^*) \leq d_0\zeta_n^*/b \leq \bar{d}_b\zeta_n^*/b$. We can easily modify Lemma A.5 to the new situation, to see that

$$2c_1(1+b)c_2\bar{\lambda}_{n,0} \leq \mathcal{V}_{\theta_n^*} + \delta\tilde{\mathcal{E}} + \delta\hat{\mathcal{E}}$$

$$= \delta\tilde{\mathcal{E}} + \epsilon_n^*.$$

So

$$\tilde{\mathcal{E}} + c_1(1+b)\bar{\lambda}_{n,0}\hat{I}(\tilde{\theta} - \theta_n^*) \leq \bar{\lambda}_{n,0}\frac{d_0\zeta_n^*}{b} + \delta\tilde{\mathcal{E}} + \epsilon_n^*$$

$$= (d_0 + b)\bar{\lambda}_{n,0}\frac{\zeta_n^*}{b} + \delta\tilde{\mathcal{E}}.$$

But then, using $\hat{I}(\tilde{\theta} - \theta_n^*) \geq I(\tilde{\theta} - \theta_n^*)/c_1$, and $0 < \delta < 1$, we derive

$$(1+b)\bar{\lambda}_{n,0}I(\tilde{\theta} - \theta_n^*) \leq c_1(d_0 + b)\bar{\lambda}_{n,0}\frac{\zeta_n^*}{b}. \qquad \square$$

Recall the definitions

$$\delta_1 := (1/(1+b))^{N_1}, \qquad \delta_2 := (1/(1+b))^{N_2}$$

and, as in (16) and (17),

$$d(\delta_1, \delta_2) := 1 + \left(\frac{1 + (d^2 - 1)\bar{\delta}_1}{(1-d)(1-\delta_1)}\right)\delta_2$$

and

$$\Delta(b, \delta, \delta_1, \delta_2) := d(\delta_1, \delta_2)\frac{1 - \delta^2}{\delta b} \vee 1.$$

**THEOREM A.5.** *Suppose Conditions* I$(b, \delta, c_1, c_2)$, II$(b, \delta, d, c_1, c_2)$ *and* III$(c_0)$ *are met, with* $c_2 = \sqrt{2c_1^2 - 1}$ *and* $c_1 > c_0 > 1$. *Take*

$$\bar{\lambda}_{n,0} > 4\sqrt{\frac{\log(2m)}{n}}\left(\sqrt{2} + \frac{\sqrt{6c_0^2 - 4} - \sqrt{2c_0^2}}{2c_0}\right).$$

*Then there exists a solution* $s > 0$ *of* $1 - 1/c_1^2 = K_m\lambda_{n,0}(s)$, *and a solution* $t > 0$ *of* $\bar{\lambda}_{n,0} = \bar{\lambda}_{n,0}(t)$. *We have with probability at least* $1 - \alpha$, *with*

$$\alpha = \exp[-na_n^2 s^2] + \left(\log_{1+b}\frac{(1+b)^2\Delta(b, \delta, \delta_1, \delta_2)}{\delta_1\delta_2}\right)\exp[-n\bar{a}_n^2 t^2],$$

*that*

$$\mathcal{E}(f_{\hat{\theta}_n}) \leq \frac{\epsilon_n^*}{1 - \delta},$$

*and moreover*

$$I(\hat{\theta}_n - \theta_n^*) \leq d(\delta_1, \delta_2)\frac{\zeta_n^*}{b}.$$

PROOF. Since $c_1 > c_0$, there exists a solution $s > 0$ of $1 - 1/c_1^2 = K_m \lambda_{n,0}(s)$. By Lemma A.9, the set $\Omega$ has probability at least $1 - \exp[-na_n^2 s^2]$. We may therefore restrict attention to the set $\Omega$.

We also know that

$$a_n = \sqrt{\frac{\log(2m)}{n}}\left(\sqrt{2} + \sqrt{\frac{\log(2m)}{n}}K_m\right)$$

$$\leq \sqrt{\frac{\log(2m)}{n}}\left(\sqrt{2} + \frac{\sqrt{6c_0^2 - 4} - \sqrt{2c_0^2}}{2c_0}\right).$$

Hence there is a solution $t > 0$ of $\bar{\lambda}_{n,0} = \lambda_{n,0}(t)$.

This means that the rest of the proof is similar to the proof of Theorem A.4, using that on $\Omega$, $I_l(\theta|\theta_n^*)/c_1 \leq \hat{I}_l(\theta|\theta_n^*) \leq c_2 I_l(\theta|\theta_n^*)$, for $l = 1, 2$. $\square$

PROOF OF THEOREM 2.2. Again, in Theorem 2.2 we have stated a special case with $b = 1$, $c_1 = 3/2$, $\delta = 1/2$, $d = 2$, $\delta_1 = \delta_2 = 1/2$. We moreover presented some simpler conservative estimates for the expressions that one gets from inserting these values in Theorem A.5. $\square$

When $K_m$ is not known, the values $s$ and $t$ in Theorem A.5 are also not known. We may however estimate them, as well as $a_n$ (and $\bar{a}_n = 4a_n$) to obtain estimates of the levels $\exp[-na_n^2 s^2]$ and $\exp[-n\bar{a}_n^2 t^2]$. Define

$$\bar{a}_n(K) := 4a_n(K), \qquad a_n(K) := \sqrt{\frac{2\log(2m)}{n}} + \frac{\log(2m)}{n}K.$$

Let

$$\hat{K}_m := \max_{1 \leq k \leq m} \frac{\|\psi_k\|_\infty}{\hat{\sigma}_k}.$$

Note that on $\Omega$,

$$K_m/c_2 \leq \hat{K}_m \leq c_1 K_m.$$

Condition III($c_1, c_2$) is somewhat stronger than Condition III($c_0$).

CONDITION III($c_1, c_2$).   For some known constant $c_1 > 1$ and for $c_2 = \sqrt{2c_1^2 - 1}/c_1$, it holds that

$$2\sqrt{\frac{\log(2m)}{n}}K_m < \frac{\sqrt{6c_1^2 - 4} - \sqrt{2c_1^2}}{c_1^2 c_2}.$$

LEMMA A.11.   *Assume Condition* III($c_1, c_2$). *On* $\Omega$, *there exists a solution* $\hat{s} \leq s$ *of the equation*

$$
\begin{aligned}
1 - \frac{1}{c_1^2} &= c_2 \hat{K}_m a_n(c_2 \hat{K}_m) \\
(20) \\
&\quad \times \left( 1 + \hat{s}\sqrt{2(1 + 2c_2\hat{K}_m a_n(c_2\hat{K}_m))} + \frac{2c_2 \hat{K}_m a_n(c_2 \hat{K}_m)\hat{s}^2}{3} \right).
\end{aligned}
$$

*Take*

$$\bar{\lambda}_{n,0} > 4\sqrt{\frac{\log(2m)}{n}}\left( \sqrt{2} + \frac{\sqrt{6c_1^2 - 4} - \sqrt{2c_1^2}}{2c_1^2 c_2} \right).$$

*Then there is a solution* $\hat{t} \leq t$ *of the equation*

$$\bar{\lambda}_{n,0} = \bar{a}_n(c_2 \hat{K}_m)\left( 1 + \hat{t}\sqrt{2(1 + 2c_2\hat{K}_m \bar{a}_n(c_2\hat{K}_m))} + \frac{2c_2 \hat{K}_m \bar{a}_n(c_2 \hat{K}_m)\hat{t}^2}{3} \right).$$

PROOF.   We have

$$2\sqrt{\frac{\log(2m)}{n}}\hat{K}_m \leq c_1 2\sqrt{\frac{\log(2m)}{n}}K_m < \frac{\sqrt{6c_1^2 - 4} - \sqrt{2c_1^2}}{c_1 c_2}.$$

Hence

$$2c_2\hat{K}_m\sqrt{\frac{\log(2m)}{n}} < \frac{\sqrt{6c_1^2 - 4} - \sqrt{2c_1^2}}{c_1}.$$

So

$$c_2 \hat{K}_m a_n(c_2 \hat{K}_m) \leq 1 - \frac{1}{c_1^2}.$$

Hence there exists a solution $0 < \hat{s} \leq s$ of (20). The solution $\hat{t}$ exists for similar reasons.   □

COROLLARY A.3.   *Suppose Condition* III($c_1, c_2$). *Let the level* $\alpha$ *be defined in Theorem* A.5 *and* $\hat{s}$ *and* $\hat{t}$ *be given in Lemma* A.11. *Define the estimated level*

$$
\begin{aligned}
\hat{\alpha}^U &= \exp[-na^2(\hat{K}_m/c_1)\hat{s}^2] \\
&\quad + \left( \log_{1+b}\frac{(1+b)^2 \bar{\Delta}(b, \delta, \bar{\delta}_1, \bar{\delta}_2)}{\bar{\delta}_1 \bar{\delta}_2} \right)\exp[-n\bar{a}_n^2(\hat{K}_m/c_1)\hat{t}^2].
\end{aligned}
$$

*Then with probability at least $1 - \exp[-na_n^2 s^2]$, it holds that $\hat{\alpha}^U \geq \alpha$.*

REMARK. Using similar arguments, one obtains an estimate $\hat{\alpha}^L$ such that $\hat{\alpha}^L \leq \alpha$ with probability at least $1 - \exp[-na_n^2 s^2]$.

### A.5. Proof of the result of Example 4.

LEMMA A.12. *Let $\mathbf{Z}(M)$ be the random variable defined in (14). Assume that $\|f_{\theta^*} - \bar{f}\|_\infty \leq \eta$ and that $MK_m + 2\eta \leq c_3$. Then we have*

$$\mathbf{P}(\{\mathbf{Z}(M) \geq \tilde{\lambda}_{n,0}M\} \cap \Omega) \leq 2\exp[-n\bar{a}_n^2 t^2],$$

*where*

$$\tilde{\lambda}_{n,0} := \tilde{\lambda}_{n,0}(t) := c_2 \sqrt{\frac{2\log(2m)}{n\bar{a}_n^2} + 2t^2} + c_3\lambda_{n,0}(t).$$

PROOF. Clearly,

$$\left| \frac{1}{n}\sum_{i=1}^n \varepsilon_i (f_{\theta^*} - f_\theta)(X_i) \right| \leq I(\theta^* - \theta) \max_{1 \leq k \leq m} \left| \frac{1}{n}\sum_{i=1}^n \varepsilon_i \psi_k(X_i) \right|.$$

Because the errors are $\mathcal{N}(0,1)$, we get for

$$a = c_2 \bar{a}_n \sqrt{\frac{2\log(2m)}{n\bar{a}_n^2} + 2t^2},$$

that

$$\mathbf{P}\left( \left\{ \left| \frac{1}{n}\sum_{i=1}^n \varepsilon_i \psi_k(X_i) \right| \geq a \right\} \cap \Omega \right) \leq 2m \exp\left[ -\frac{na^2}{2c_2^2} \right] = \exp[-n\bar{a}_n^2 t^2].$$

Moreover, the function $x \mapsto x^2/(2c_3)$ is Lipschitz when $|x| \leq 1$. Since $\|f_\theta + \bar{f}\|_\infty \leq 2\eta + MK_m \leq c_3$, we can apply Corollary A.1 to find that with probability at least $1 - \exp[-n\bar{a}_n^2 t^2]$, we have

$$\sup_{f \in \mathcal{F}_M} \tfrac{1}{2}|(Q_n - Q)((f_\theta - \bar{f})^2 - (f_{\theta^*} - \bar{f})^2)| \leq c_3 M\lambda_{n,0}(t). \qquad \square$$

PROOF OF THEOREM 3.1. By Lemma A.9,

$$\mathbf{P}(\Omega) \geq 1 - \exp[-na_n^2 s^2].$$

Using Lemma A.12 and applying the same arguments as in the proof of Theorem A.4, the result follows for general constants $b$, $\delta$, $d$, $\delta_1$ and $\delta_2$ and constants $c_1$, $c_2 = \sqrt{2c_1^2 - 1}/c_1$ and $c_3$. Theorem 3.1 takes $b = 1$, $\delta = 1/2$, $d = 2$, $\delta_1 = \delta_2 = 1/2$ and $c_1 = 3/2$, $c_2 := \sqrt{14/9}$ and $c_3 = 1$. $\quad\square$

SEMINAR FÜR STATISTIK
ETH ZÜRICH
LEO D11
8092 ZÜRICH
SWITZERLAND
E-MAIL: geer@stat.math.ethz.ch